\documentclass [11pt]{article}
\usepackage{latexsym,amsfonts}
\begin{document}
\bibliographystyle{plain}

\newtheorem{defi}{Definition}[section] 
\newtheorem{theo}{Theorem}[section]
\newtheorem{prop}{Proposition}[section] 
\newtheorem{lemm}{Lemma}[section]     
\newtheorem{conj}{Conjecture}[section] 
\newtheorem{note}{Note}[section] 
\newtheorem{coro}{Corollary}[section]    
\newtheorem{ques}{Question}[section] 
\newtheorem{rema}{Remark}[section]

\newcommand{\fas}{\ensuremath{f_{\rm {A, sol}}(n)}} 
\newcommand{\ms}{\ensuremath{M_{\rm{s,p',A}}(G)}} 
\newcommand{\fbh}{\ensuremath{f({\mathbb Z}_{p^\beta}H, \alpha)}}   
\newcommand{\fah}{\ensuremath{f({\mathbb Z}_{p^\alpha}H, \alpha)}}   
\newcommand{\zah}{\ensuremath{{\mathbb Z}_{p^\alpha}H}} 
\newcommand{\zbh}{\ensuremath{{\mathbb Z}_{p^\beta}H}} 
\newcommand{\zih}{\ensuremath{{\mathbb Z}_{p^\alpha}H_i}} 
\newcommand{\zpq}{{\mathbb Z}_pQ} 
\newcommand{\zph}{{\mathbb Z}_pH} 
\newcommand{\phx}{\ensuremath{{\phi}_{\scriptscriptstyle X}}}   
\newcommand{\jpq}
{\ensuremath{{\mathfrak A}_p{\mathfrak A}_q \vee {\mathfrak A}_q{\mathfrak A}_p}} 
\newcommand{\pqr}
{\ensuremath{{\mathfrak A}_p{\mathfrak A}_q{\mathfrak A}_r}}
\newcommand{\jpqr}
{\ensuremath{ \bigvee_{\pi}\left\{{\mathfrak A}_{\pi(p)}{\mathfrak A}_{\pi(q)}
{\mathfrak A}_{\pi(r)}\right\}}}
\newcommand{\fu}{\ensuremath{f_{\mathfrak U}(n)}}
\newcommand{\fv}{\ensuremath{f_{\mathfrak V}(n)}}
\newcommand{\fw}{\ensuremath{f_{\mathfrak W}(n)}} 
\newcommand{\fvp}{\ensuremath{f_{\mathfrak V}(p^{\alpha}q^{\beta}r^{\gamma})}} 
\newcommand{\fxp}{\ensuremath{f_{\cal X}(p^{{\alpha}_1}q^{{\beta}_1}r^{{\gamma}_1})}} 
\newcommand{\fy}{\ensuremath{f_{\cal Y}(p^{\alpha}_2q^{\beta}_2)}} 
\newcommand{\fbx}{\ensuremath{f_{\bar {\cal X}}(p^{\alpha}_1q^{\beta}_1)}} 
\newcommand{\fby}{\ensuremath{f_{\bar {\cal Y}}(p^{\alpha_2}q^{\beta_2})}} 

\newcommand{\lav}{\ensuremath{a = \limsup\left \{\frac{\log f_{\mathfrak V}(n)}{\mu (n) \log n}\right \} }}
\newcommand{\lau}{\ensuremath{a = \limsup\left \{\frac{\log f_{\mathfrak U}(n)}{\mu (n) \log n}\right \} 
}}

\newcommand{\law}{\ensuremath{a = \limsup\left \{\frac{\log f_{\mathfrak W}(n)}{\mu (n) \log n}\right \} 
}}

\title{Groups in which Squares and Cubes Commute} 

\author{Geetha Venkataraman \\
School of Liberal Studies \\
Ambedkar University Delhi\\
Lothian Road\\
Kashmere Gate\\
Delhi-110006 \\
India \\
\vspace{.05in}
email: geetha@aud.ac.in}
\date{\today}      
\maketitle

\noindent {\bf Abstract}
If for all $a, b$ in a group $G$,  we have that $a^2b^2 = b^2a^2$ and $a^3b^3 = b^3a^3$ then does the group necessarily have to be abelian? This paper shows that the answer is affirmative for finite groups  as well as certain classes of infinite groups. In general we show that if $G$ is an infinite group in which squares commute and cubes commute then the set of torsion elements of $G$ denoted by $T(G)$ has to be an abelian normal subgroup of $G$.

\vspace{.15in}

\vspace{.15in}

\section{Introduction}  

There are various conditions under which a group becomes abelian. It is a trivial exercise to show that in a group $G$ if ${(ab)}^2 = a^2b^2$ for all $a, b \in G$ then $G$ is abelian. Another condition which also results in the group being abelian is if for all $a, b \in G$ we have ${(ab)}^r = a^r b^r$ for three consecutive integers $r$.  

It is natural to wonder about a group $G$ in which $a^2b^2$ = $b^2a^2$ for all $a, b \in G$. However this is not sufficient to ensure that $G$ is abelian. An example of such a group is $S_3$, the group of permutations on $3$ letters. The squares of elements in $S_3$ form a cyclic group generated by a $3$-cycle and so it is clear that squares of elements commute in $S_3$. However $S_3$ is not abelian. As a next step, one can ask what happens if in addition to squares commuting, cubes also commute.

We are able to prove using results from the theory of finite groups and Sylow subgroups that any finite group in which $a^2b^2 = b^2a^2$ and $a^3b^3 = b^3a^3$ has to be abelian. This is one of the main results of the paper. 

It is also shown that the numbers $2, 3$ are not important and that there is a more general result. If $m,n$ are coprime natural numbers and if $a^mb^m = b^ma^m$ and $a^nb^n = b^n a^n$ for all $a, b$ in a finite group $G$, then $G$ must be abelian.

Some definitions and terminology that we shall use in the paper are given below.

Let $m, n$ be fixed coprime natural numbers and let us say that a group $G$ has property $\cal P$ if  $a^ib^i = b^ia^i$ for $i=m, n$ and for all $a, b \in G$. Note that if $G$ satisfies $\cal P$, then so do subgroups and quotients.

Let ${\cal X}$ be a class of groups. We say that $G$ is a $\cal X$-group if $G \in \cal X$. Some definitions: a group $G$ is {\bf locally $\cal{X}$} if every finitely generated subgroup of $G$ is in ${\cal X}$ and a group $G$ is {\bf residually $\cal X$} if for each non-identity element $g \in G$  there is a normal subgroup $N_g$ of $G$ such that $g \not \in N$ and ${G/ N_g \in \cal X}$.


The main results we prove are given below.

\vspace{.05in}
\noindent{\bf Theorem 2.1}
{\em Let ${\cal X} $ be a class of groups. Suppose that all $\cal{X}$-groups satisfying $\cal P$ are abelian.  Then the following statements hold.}
\begin{enumerate}
\item[(a)] {\em All locally $\cal X$-groups satisfying $\cal P$ are abelian.}
\item[(b)] {\em All residually $\cal X$-groups satisfying $\cal P$ are abelian.}
\end{enumerate}
\vspace{.05in}

\vspace{.05in}
\noindent{\bf Theorem 3.1}
{\em Let $G$ be a finite group satisfying $\cal P$. Then $G$ is abelian.}
\vspace{.05in}

\vspace{.05in}
\noindent{\bf Theorem 4.1}
{\em Let $G$ be an infinite group satisfying $\cal P$.} 
\begin{enumerate}
\item[(a)] {\em Then $T(G)$ the set of all torsion elements of $G$ is a normal abelian subgroup of $G$.}
\item[(b)] {\em Any Sylow subgroup of $G$ is abelian and in particular, if $G$ is an infinite $p$-group then $G$ is abelian.}
\item[(c)] {\em If $G$ is locally finite or residually finite then $G$ is abelian.}
\item[(d)] {\em If $G$ is residually $p$ for a prime $p$ then $G$ is abelian.}

\end{enumerate}
\vspace{.05in}

The rest of the paper is organised as follows. The next section has the preliminary results required to prove the main results and the proof of Theorem 2.1. The third section has the proof for the finite case. The fourth section deals essentially with proving the main result concerning infinite groups. The last section has concluding remarks and discusses some open questions.

For proofs of well known theorems that are used here and terminology see \cite{djsr95} or \cite{jjr94}.

\section{Preliminary Results}

\begin{lemm}
Let $G$ be a group satisfying $\cal P$.  Let $k$ be a natural number that is coprime to at least one of $m$ or $n$. If $H = \{ x \in G \mid x^k =e\}$ then $H$ is a normal abelian subgroup of $G$.
\end{lemm}

\noindent{\bf Proof}{   } It is clear that if $H$ is a subgroup of $G$ then it will be a normal subgroup since for any $x \in H$ and $g \in G$, ${(gxg^{-1})}^k = gx^kg^{-1} =e$. 

Note that $e \in H$ and that $H$ is closed under taking inverses. Thus we just have to show closure and we do so by first showing that elements of $H$ commute. Consider any $x, y \in H$. Without loss of generality let us assume that $k$ is coprime to $m$. Thus there exist integers $\lambda$ and $\mu$ such that $\lambda k + \mu m = 1$. So for any $h \in H$ we have $h = {(h^k)}^\lambda {(h^\mu)}^m = {(h^\mu)}^m$. So for any $x, y \in H$, using property $\cal P$, we have
\begin{eqnarray*}
xy & =& {(x^\mu)}^m {(y^\mu)}^m \\
     & = & {(y^\mu)}^m{(x^\mu)}^m \\
     & = & yx.
     \end{eqnarray*}   
Consequently for any $x, y \in H$, we have ${(xy)}^k = x^ky^k = e$. Thus $H$ is an abelian normal subgroup of $G$. \ $\Box$

\begin{theo} Let ${\cal X} $ be a class of groups. Suppose that all $\cal{X}$-groups satisfying $\cal P$ are abelian.  Then the following statements hold.
\begin{enumerate}
\item[{\rm(a)}] All locally $\cal X$-groups satisfying $\cal P$ are abelian.
\item[{\rm(b)}]  All residually $\cal X$-groups satisfying $\cal P$ are abelian.
\end{enumerate}
\end{theo}

\noindent{\bf Proof}{   } Let ${\cal X} $ be a class of groups and suppose that all $\cal{X}$-groups satisfying $\cal P$ are abelian. 

Let $G$ be a locally $\cal X$-group which satisfies $\cal P$. Thus $G$ satisfies $\cal P$ and every finitely generated subgroup of $G$ is in $\cal X$. Let $x, y \in G$ and consider $ H = \langle x, y \rangle$. Then $H$ is a finitely generated subgroup of $G$ and hence must be in $\cal X$ and clearly $H$ satisfies $\cal P$. Thus $H$ is abelian and so $xy =yx$. Thus $G$ is abelian.

Now let $K$ be a residually $\cal X$-group which satisfies $\cal P$. Thus $K$ satisfies $\cal P$ and for each non-identity element $x \in K$  there is a normal subgroup $N_x$ of $K$ such that $x \not \in N$ and ${K/ N_x\in \cal X}$. Now $K/ N_x$ is a member of $\cal X$ and satisfies $\cal P$. Thus $K/N_x$ is abelian.

Let ${\hat K} = \mbox{Cr}_{e \not = x \in K} K/N_x$ denote the group which is the cartesian product or unrestricted product of the groups $K/ N_x$ as $x$ varies over $K \setminus \{e\}$. Note that ${\hat K}$ is an abelian group. We can define a map $\iota: K \rightarrow {\hat K}$ as $\iota(y) = (yN_x)$. Then $\iota$ is a monomorphism and so $K$ is a subdirect product of a cartesian product of abelian groups and hence will be abelian. \ $\Box$

\section{Finite groups satisfying ${\cal P}$}

We begin by proving the following lemma.

\begin{lemm}

Let $G$ be a finite group and let $p$ be a prime such that $p^{\alpha}$ is the highest power of $p$ dividing $|G|$. Define $G_{p^{\alpha}} = \{x \in G \mid x^{p^{\alpha}} =e\}$. Then $G_{p^{\alpha}}$ is a subgroup of $G$ if and only if $G$ has a unique  Sylow $p$-subgroup, $P$. Further, in either case, $P = G_{p^{\alpha}}$.
\end{lemm}

\noindent{\bf Proof}{   } Let $G$ be a finite group and $p$ a prime such that $p^{\alpha}$ is the highest power of $p$ dividing $|G|$ and where $\alpha \in {\mathbb N}$. Let $G_{p^{\alpha}} = \{x \in G \mid x^{p^{\alpha}} =e\}$ be a subgroup of $G$. By the defining property of $G_{p^{\alpha}}$, it is clear that $G_{p^{\alpha}}$ must have order that is a power of $p$. For if $q$ is a different prime that divides $|G_{p^{\alpha}}|$ then by Cauchy's Theorem $G_{p^{\alpha}}$ will have an element $g$ of order $q$. Also $g^{p^{\alpha}} =e$. Thus $q|p^{\alpha}$ and we get that $q=p$, a contradiction. Thus $|G_{p^{\alpha}}| = p^{\beta}$ for some $\alpha \geq \beta \geq 0$. Since $p$ divides $|G|$, by Sylow's Theorem we know that a Sylow $p$-subgroup exist. Let $P$ be any Sylow $p$-subgroup of $G$. Then $|P|= p^{\alpha}$. Thus for all $x \in P$,  we will have $x^{p^{\alpha}} = e$ and so $P \subseteq G_{p^{\alpha}}$. Thus $\alpha \leq \beta$ and so we get that $G_{p^{\alpha}} = P$. This shows that there is a unique and hence a normal Sylow $p$-subgroup in $G$, namely, $ P = G_{p^{\alpha}}$.

Conversely let $G$ have a unique Sylow $p$-subgroup, say $P$ with $|P| = p^{\alpha}$. Let $G_{p^{\alpha}} = \{ x \in G \mid x^{p^{\alpha}} =e\}$. Then as seen above we must have $P \subseteq G_{p^{\alpha}}$. Now if $x \in G_{p^{\alpha}}$ then $x^{p^{\alpha}} =e$. So $x$ has $p$ power order. Thus $\langle x \rangle$ is a $p$-subgroup of $G$ and hence must be contained in a Sylow $p$-subgroup of $G$. Since $P$ is the unique Sylow $p$-subgroup of $G$, therefore $x \in P$ and so $G_{p^{\alpha}} \subseteq P$. Thus $G_{p^{\alpha}} = P$ and so is a subgroup of $G$. \ $\Box$

Now we are in a position to present our main result concerning finite groups.

\begin{theo}
Let $G$ be a finite group satisfying $\cal P$. Then $G$ is abelian.
\end{theo}

\noindent{\bf Proof}{   } Let $G$ be a finite group satisfying $\cal P$ and let $|G| = {p_1}^{\alpha_1} \cdots {p_r}^{\alpha_r}$ be the prime decomposition of $|G|$. For each $i = 1, \ldots , r$, let $P_i$ denote the Sylow $p_i$-subgroups of $G$.  Now $G$ satisfies $\cal P$, and for each $i$, the prime $p_i$ will be coprime to either $m$ or $n$. Thus by Lemma 2.1, we will have that $G_{{p_i}^{\alpha_i}} = \{ x \in G \mid x^{{p_i}^{\alpha_i}}=e\}$ is a normal abelian subgroup of $G$. Hence by Lemma 3.1, $P_i$ is the unique normal Sylow $p_i$-subgroup of $G$. Thus $G = P_1 \times \cdots \times P_r$ as a direct product and hence abelian.\ $\Box$

\section{Infinite Groups satisfying $\cal P$}

\begin{theo}
Let $G$ be an infinite group satisfying $\cal P$.
\begin{enumerate}
\item[{\rm(a)}]  Then $T(G)$ the set of all torsion elements of $G$ is a normal abelian subgroup of $G$.
\item[{\rm (b)}]  Any Sylow subgroup of $G$ is abelian and in particular if $G$ is an infinite $p$-group then $G$ is abelian.
\item[{\rm(c)}]  If $G$ is locally finite or residually finite then $G$ is abelian.
\item[{\rm(d)}]  If $G$ is residually $p$ for a prime $p$ then $G$ is abelian.
\end{enumerate}
\end{theo}

\noindent{\bf Proof}{   }Let $G$ be an infinite group satisfying $\cal P$.

\vspace{.05in}
\noindent{(a)}{   } Let $H = T(G)$ be the set of all torsion elements of $G$.  To show that $H$ is a subgroup we only need to show closure. If $H = \{e\} $ then the result is trivially true. So assume otherwise and let $x \in H \setminus \{e\}$ and let $o(x) =  r = p_1^{\alpha_1} \cdots p_k^{\alpha_k}$ be the prime decomposition for $r$. For each $i = 1, \ldots, k$, let $$q_i = \frac{r}{p_i^{\alpha_i}}.$$ 
Then since $\gcd(q_1, \ldots, q_k) =1$ we have  $\lambda_1q_1 + \cdots + \lambda_kq_k =1$ for some $\lambda_i \in {\mathbb Z}$ such tha. Thus we can write
\begin{eqnarray*} 
x & = & x^{\lambda_1q_1 + \cdots + \lambda_kq_k } \\
  & = & x^{\lambda_1 q_1} \cdots x^{\lambda_k q_k} \\
  & = & x_1^{\lambda_1} \cdots  x_k^{\lambda_k}
  \end{eqnarray*}
where for each $i$, we define $x_i = x^{q_i}$. Then  $o(x_i) = p_i^{\alpha_i}$ and so $x_i \in H$.  Now given any $x, y \in H \setminus \{e\}$ we can write $x = x_1^{\lambda_1} \cdots x_k^ {\lambda_k}$ and $y = y_1^{\mu_1} \cdots y_t^{\mu_t}$ where $x_i$ and $y_j$ are elements of $H$ of prime power orders for all $i, j$. So if we show that elements of prime power order in $H$ commute then elements of $H$ will commute. Thus closure will hold in $H$. Consequently $H$ will be an abelian subgroup of $G$. It is obvious that if $T(G) \leq G$ then $T(G)$ is normal. 

Now consider any $a, b \in H$ where $a$ and $b$ have prime power order say $p^\alpha$ and $q^\beta$ respectively. Let $s = p^\alpha q^\beta$. If $s$ is coprime to at least one of $m$ or $n$, then by Lemma 2.1,  $S = \{ g \in G \mid g^s =1\}$ is a (normal) abelian subgroup of $G$ and clearly $a, b  \in S$. Thus $ab =ba$. 

So assume that $s$ is neither coprime to $m$ nor to $n$. Then since $\gcd(m, n) = 1$ we must have that $p \not = q$ and without loss of generality let us assume that $p |m$, $q|n$. Now let $S_a = \{ g \in G \mid g^{p^\alpha} =1\}$ and $S_b = \{ g \in G \mid g^{q^\beta} =1\}$. Then since $p^\alpha$ is coprime to $n$ and $q^\beta$ is coprime to $n$, by Lemma 2.1 we get that $S_a$ and $S_b$ are both normal abelian subgroups of $G$. Further $a \in S_a, b \in S_b$. Also since $p$ and $q$ are distinct primes $S_a \cap S_b = \{e\}$. But then normality of the subgroups ensures that $gh = hg$ for all $ g \in S_a$ and $h \in S_b$. Consequently $ab=ba$.

Since elements of prime power in $H$ commute we can conclude that $ T(G) = H$ is an abelian (normal) subgroup of $G$.

\vspace{.05in}
\noindent{(b)}{   }This is just a special case of (a), for if $P$ is a Sylow $p$-subgroup of $G$ for some prime $p$, then every element of $P$ has $p$-power order and hence finite order. It is equally clear for the case when $G$ is an infinite $p$-group.

\vspace{.05in}
\noindent{(c)}{   }Let ${\cal X}$ be the class of finite groups. Then by Theorem 3.1, all ${\cal X}$ groups satisfying $\cal P$ are abelian. Therefore by Theorem 2.1 all locally $\cal X$-groups and residually $\cal X$-groups satisfying $\cal P$ are abelian. Hence if $G$ is a locally finite group or a residually finite group then $G$ is abelian.

\vspace{.05in}
\noindent{(d)}{  }Let $p$ be a prime and let ${\cal X}$ be the class of $p$-groups. Then by Theorem 3.1 and part (b) above, all ${\cal X}$ groups satisfying $\cal P$ are abelian. Therefore by Theorem 2.1 all residually $\cal X$-groups satisfying $\cal P$ are abelian. Hence if $G$ is a residually $p$-group then $G$ is abelian. \ $\Box$.

\begin{rema}
{\rm
\begin{enumerate}
\item [(i)] The proof given for Theorem 4.1 (a) is another proof for the result that every finite group $H$ satisfying $\cal P$ is abelian. Indeed the proof also shows that if $G$ is any group satisfying $\cal P$ and $T$ is any subset of torsion elements of $G$ that contains identity and is closed with respect to inverses then it must be an abelian subgroup of $G$.
\item [(ii)] It can be shown independently that if $G$ is locally nilpotent then $T(G) \leq G$. Clearly by Theorem 4.1 (a), if such a $G$ satisfies $\cal P$ then $T(G)$ will be an abelian normal subgroup of $G$. However we include below a proof which discusses the special case of locally nilpotent groups. 

Let $G$ be locally nilpotent group satisfying $\cal P$. Let $T(G)$ denote the set of torsion elements of $G$. Let $x, y \in T(G)$ and let $H = \langle x, y \rangle$. Note that $T(G)$ contains $e$ and is closed under taking inverses. Indeed if $T(G) \leq G$ then it has to be normal.

Since $G$ is locally nilpotent therefore $H$ is nilpotent and is generated by two elements of finite order. Hence $H$ is finite and so $xy \in T(G)$. So $T(G) \leq G$. Now $H$ is a finite group satisfying $\cal P$. Thus by Theorem 3.1, $H$ is abelian. So $xy =yx$. Therefore  $T(G)$ is also abelian.

\item [(iii)] Since every nilpotent group is locally nilpotent, if $G$ is a nilpotent group satisfying $\cal P$, then by the above remark (and also Theorem 4.1 (a)), $T(G)$ is an abelian normal subgroup of $G$.
\item [(iv)] Since every finitely generated torsion-free nilpotent group $G$ is a residually finite $p$-group for every prime $p$ by Theorem 4.1 (c) if $G$ also satisfies $\cal P$ then $G$ will be abelian. 
\item [(v)] If $G$ is a polycyclic group then it is residually finite and so if $G$ satisfies $\cal P$ then by Theorem 4.1 (c), $G$ is abelian.
\end{enumerate}}
\end{rema}

\section{Conclusion}

From the earlier sections, we can see that a finite group satisfying $\cal P$ is abelian. However the question as to whether an infinite group satisying $\cal P$ is abelian still remains open. 

It is clear from Theorem 4.1 (a) that if $G$ is an infinite group whose torsion elements do not form a group or form a non-abelian group then such a group cannot satisfy $\cal P$ for any choices of coprime $m, n$.

Thus far, however we have neither been able prove that an infinite group $G$ satisfying $\cal P$ is abelian nor provide a counter example. The counter example would have to be a non-abelian infinite group with an abelian torsion subgroup and which satisfies $\cal P$ for some coprime natural numbers $m, n$. One could also try to see if there are torsion free non-abelian infinite groups that satisfy $\cal P$ for some coprime natural numbers $m, n$. 
\bibliography{vision1}
\end{document}